\tolerance=20000
\magnification=\magstep1
\mathsurround=2pt
\footline={\ifnum\pageno=1{\hfil}\else{\hss\tenrm\folio\hss}\fi}
\font\bigbf=cmbx10 scaled \magstep1
\font\bigit=cmti10 scaled \magstep1
\font\biggtenex=cmex10 scaled \magstep2
\font\tenmsy=msym10
\font\sevenmsy=msym7
\font\fivemsy=msym5
\newfam\msyfam
\textfont\msyfam=\tenmsy
\scriptfont\msyfam=\sevenmsy
\scriptscriptfont\msyfam=\fivemsy
\def\Bbb#1{\fam\msyfam#1}
\baselineskip=16pt
\lineskip=3pt
\lineskiplimit=3pt
\def\cW{{\cal W}}
\def\R{{\Bbb R}}
\def\topsim#1{\ \smash{\mathop{\sim}\limits^{#1}}\ }
\def\ttopsim#1{\ {\buildrel{#1}\over\sim}\ }
\def\bigk{\mathop{\biggtenex K}}

\def\square{\vbox{
      \hrule height 0.4pt
      \hbox{\vrule width 0.4pt height 5.5pt \kern 5.5pt \vrule width 0.4pt}
      \hrule height 0.4pt}}

\def\Pn1{P_{n}^{(1)} (z)}
\def\Pn0{P_{n}^{(0)} (z)}

\def\tphin{\> {}_{10}\phi_9}
\def\ephis{\> {}_8\phi_7}

\def\fphit{\> {}_4\phi_3}
\def\sFs{\> {}_7 F_6 }

\def\bibline{\ {\vrule width.5in height .03pt depth .05pt}\ }
\rm
\vglue .5in
\centerline {\bigbf Solutions to the}
\centerline{\bigbf Associated  {\bigit q }-Askey-Wilson
Polynomial Recurrence Relation \footnote{*}{{\rm Research partially
supported by NSERC (Canada).}}}
\vglue .25in
\centerline {Dharma P. Gupta and David R. Masson}
\bigskip
\baselineskip=12pt
\centerline {Department of Mathematics}
\centerline {University of Toronto}
\centerline {Toronto Ontario M5S 1A1}
\centerline {Canada}
\bigskip
\centerline{{\it In honour of Walter Gautschi on the occasion of 
his 65th birthday.}}
\vglue .75in

\baselineskip=14pt
\midinsert\narrower\narrower
\centerline {\bf Abstract}
\medskip

A  $ \tphin $  contiguous relation is used to derive contiguous relations
for a very-well-poised  $ \ephis $.  These in turn yield solutions
to the associated  $ q $-Askey-Wilson polynomial recurrence relation,
expressions for the associated continued fraction, the weight
function and a  $ q $-analogue of a generalized Dougall's theorem.
\bigskip

\medskip

\endinsert
\vfill\eject

\baselineskip=18pt
\beginsection  1. Introduction.

A deeper insight into the properties of classical orthogonal
polynomials is obtained by studying their associated cases.
For recent work on associated
Hermite, Laguerre, ultraspherical,
Jacobi, continuous Hahn, continuous dual Hahn,
big  $ q $-Laguerre and big
$ q $-Jacobi polynomials see [2], [4], [8], [9], [11], [12], [14], [15],
[23].  The most general case of associated  $ q $-Askey-Wilson
polynomials [16] and their  $ q \to 1 $  limit, the associated
Wilson polynomials, have also investigated  [13], [18].

In this paper we examine the associated $ q $-Askey-Wilson 
polynomial case using the methods of [18]. All of our parameters
are, in general, complex. We essentially follow the notation in Gasper and
Rahman [5], except that we omit the designation $ q $ for the base
in the q-shifted factorials and basic hypergeometric functions.
Thus we have, with $ | q | < 1 $,
$$
(a)_0:=1, \quad (a)_n:=\prod_{j=1}^n(1-aq^{j-1}),\quad n>0, \hbox
{ or }  n=\infty,
$$
$$
(a_1, a_2, \cdots, a_k)_n:= \prod_{j=1}^k (a_j)_n .
$$
The basic hypergeometric series $ _{r+1} \phi_r $ is
$$
_{r+1} \phi_r\left( \matrix{ a_1,&a_2,&\cdots,&a_{r+1} \cr
b_1,& b_2,&\cdots
,&b_r \cr}
\quad;z\right):=\sum_{n=0}^{\infty} {(a_1,a_2, \cdots,a
_{r+1})_n \over(b_1,b_2,\cdots,b_r,q)_n}\,z^n\quad.
$$

The explicit form of monic  $ q $-Askey-Wilson polynomials is
given by [1]:
$$
\leqalignno{
P_n (z) \ = \ 
& P_n(z;\alpha,\beta,\gamma,\delta) \ = \ {\big(\alpha \beta, \alpha
\gamma, \alpha \delta \big)_n\big(\alpha \beta \gamma \delta /q\big)_n\over
(2\alpha)^n \big(\alpha \beta \gamma \delta /q\big)_{2n} 
} & (1.1)\cr
\noalign{\medskip}
& \times \fphit
\left( {{q^{-n}, \alpha \beta \gamma \delta q^{n-1}, \alpha u, \alpha /u}\atop
{\alpha \beta, \alpha \gamma, \alpha \delta}}; \ q \right) , \qquad z \ = \ {u + u^{-1}\over 2} \cr}
$$
which satisfy the three-term recurrence
$$
\leqalignno{
& P_{n+1} (z) - (z-a_n) P_n (z) + b_n^2 P_{n-1} (z) \ = \ 0 \ ,
 & (1.2) \cr
& a_n \ = \  a_n (\alpha, \beta, \gamma, \delta )
\ = \ -A_n - B_n + {\alpha
\over 2} + {1\over 2\alpha} \ , \cr
& b_n^2 \ = \  b_n^2 (\alpha, \beta, \gamma, \delta ) 
\ = \ A_{n-1} B_n \ , \cr
& A_n \ = \ { \big(1 -\alpha \beta \gamma \delta q^{n-1}\big) 
( 1- \alpha \beta q^n )
( 1- \alpha \gamma q^n )
( 1- \alpha \delta q^n)\over
2\alpha (1-\alpha \beta \gamma \delta q^{2n-1}) 
(1- \alpha \beta \gamma \delta q^{2n}) } \ , \cr
& B_n \ = \ { \alpha (1-q^n) (1 -\beta \gamma q^{n-1}) 
( 1- \beta \delta q^{n-1})
( 1- \gamma \delta q^{n-1} )\over
2 (1-\alpha \beta \gamma \delta q^{2n-2}) 
(1- \alpha \beta \gamma \delta q^{2n-1}) } \ . \cr}
$$

Associated monic  $ q $-Askey-Wilson polynomials are defined [16]
by introducing an extra parameter, say  $ \epsilon $, through a 
shift the discrete variable  $ n $  in the coefficients of (1.2).
Here we choose to replace  $ n $  by  $ n  + (\log \epsilon)/(\log q) $
so that $ q^n $ is replaced by $ \epsilon q^n $.
The new polynomial set
$ \big\{ P_n (z; \epsilon)\big\} $  is an initial value solution
$ (X_{-1} = P_{-1} = 0 , \ X_0 = P_0 = 1 ) $  to the
difference equation
$$
X_{n+1} - (z - a_v) X_n + b_v^2 X_{n-1} \ = \ 0 \qquad
v := n + (\log \epsilon)/(\log q) .
\leqno(1.3)
$$

The associated  $ q $-Askey-Wilson polynomials have been studied
by Ismail and Rahman [16].  They derive contiguous relations for a
very well poised  $ \ephis $ series and then obtain two linearly
independent solutions of the three-term recurrence relation for the 
associated Askey-Wilson polynomials.  Explicit representations of
two families of associated Askey-Wilson polynomials are given. 
Corresponding continued fractions are identified and weight 
functions obtained.

On the other hand Masson [18] studied the associated Wilson
polynomials by obtaining  $ \sFs $  contiguous relations with a
procedure (see Wilson [22]) different from that employed by
Ismail and Rahman.  He discussed explicit forms of solutions
to the three-term recurrence satisfied by associated Wilson
polynomials and identified minimal solutions for different regions.
Using Pincherle's theorem [6],[17], 
he obtained a number of new continued
fraction representations and discussed the spectral properties of the
corresponding Jacobi matrix.  A weight function was obtained and his
method gave a generalization of Dougall's theorem [3] which removed both
the terminating and the balancing conditions.

The object of the present paper is to study the associated  
$ q $-Askey-Wilson polynomials using Masson's approach [18].  In
spite of the fact that some of the basic results here will
overlap with the results obtained by Ismail and Rahman [16], we
have to offer a different procedure for obtaining the $ \ephis $  
contiguous relations, the minimal solution to the associated Askey
-Wilson poloynomial recurrence, the corresponding continued fraction
representations and an extension of Masson's generalization of
Dougall's theorem [18] to the basic hypergeometric case.  Moreover
we again see, as in many previous cases ([7], [8], [9], [18], [19],
[21]),
that Pincherle's theorem [6],[17] plays a crucial role. The importance
of Pincherle's theorem was re-emphasized by Gautschi in 1967 [6].

In Section 2, following Wilson's method [22], we derive a pair of
three-term recurrence relations for very-well-poised balanced and
terminating
$ \tphin $  basic hypergeometric functions.  From these recurrence
relations and using certain known transformations [5] we obtain a
number of explicit solutions to (1.3).  Minimal solutions for
different regions of the complex plane are identified.  In Section 3,
the associated continued fractions are obtained. In Section 4 we 
derive an
expression for the weight function of orthogonality and 
a Wronskian identity which turns out to be the basic analogue of
a generalized Dougall's theorem [18].
\vglue .5truein

\beginsection 2. Recurrence relation and solutions.

In a recent paper [10] we have derived  contiguous relations for a very- 
well-poised balanced  $ \tphin $  basic hypergeometric function
$$
\phi\ := \  
\tphin \left({ {a, q\sqrt a , -q \sqrt a , b, c, d, e, f, g, h  }
\atop {\displaystyle\sqrt a, -\sqrt a, {aq\over b}, {aq\over c}, \ldots ,
{aq\over h}}} ; \ q\right) 
\leqno(2.1)
$$
where one of the numerator parameters, say  $ h $, is equal to
$ q^{-n} $,  $ n = 0,1,\ldots $,  $ bcdefgh = a^3 q^2 $  and
$ |q| < 1 $.  One such relation is
$$
\leqalignno{
& { b(1-c) (1-{a\over c}) 
(1- {aq\over c})(1-{aq\over bd}) 
(1- {aq\over be})
(1- {aq\over bf})
(1- {aq\over bg})
(1- {aq\over bh})\over
(1- {cq\over b}) } & (2.2)\cr
&\hskip  2in \times \big[\phi (b -, c+) -\phi \big] \cr
- & { c(1-b) (1-{a\over b}) 
(1- {aq\over b}) 
(1- {aq\over cd})
(1- {aq\over ce})
(1- {aq\over cf})
(1- {aq\over cg})
(1- {aq\over ch})\over
(1- {bq\over c}) } \cr
&\hskip  2in \times \big[\phi (c -, b+) -\phi \big] \cr
- & {aq\over c} (1-{c\over b}) (1-{bc\over aq})
(1-d) (1-e) (1-f) (1-g) (1-h) \phi \ = \ 0 \ .\cr}
$$
Note that we are using the notation $\phi (b -, c +)$ to mean
$\phi$ in (2.1) with $b$ replaced by $b/q$ and $c$ replaced by
$cq$.

\noindent In (2.2), if we write  $ h = q^{-n}$, 
$ g = a^3 q^{n+2}/ bcdef $  and 
let  $ n \to \infty $,  we obtain
$$
\leqalignno{
& {bcdef\over a^2 q} 
{(1-c)  (1-{a\over c})
(1-{aq\over c})
(1-{aq\over bd})
(1-{aq\over be})
(1-{aq\over bf})\over 
(1- {cq\over b}) } \ 
\big[ \ephis (b-, c+) - \ephis \big] & (2.3) \cr-
& {bcdef\over a^2 q} 
{(1-b)  (1-{a\over b})
(1-{aq\over b})
(1-{aq\over cd})
(1-{aq\over ce})
(1-{aq\over cf})\over 
(1- {bq\over c}) } \ 
\big[ \ephis (c-, b+) - \ephis \big]  \cr
& \hskip .5truein -
 {aq\over c} (1-{c\over b}) (1-{bc\over aq})
(1-d) (1-e) (1-f) \ephis \ = \ 0 \ ,\cr}
$$
where
$$
\leqalignno{
\ephis \ = \ & \ephis 
\left (  
 { {a,q\sqrt a, -q\sqrt a, b, c,d,e,f} 
   \atop 
 {\displaystyle \sqrt a, -\sqrt a, {aq\over b}, {aq\over c},
 {aq\over d}, {aq\over e},
 {aq\over f}}} \ ; \  
  {a^2 q^2\over bcdef} \right)  \cr
  = \ & W(a; \ b, c, d, e, f ) \quad \hbox{say} \ . \cr}
$$
For non-terminating  $ \ephis $,  the convergence condition is
$ | a^2 q^2/ bcdef | < 1 $.
In (2.3) we now put
$$
\leqalignno{
a \ = \ &  {\alpha \beta \gamma u\over q}, \quad
b \ = \ q^{-n}/ \epsilon\ , \quad
c \ = \ \epsilon s q^{n-1}  \ , \quad
d \ = \ \alpha u \ , \quad
e \ = \ \beta u \ , \quad
f \ = \ \gamma u ,\cr
s \ = \ & \alpha\beta\gamma \delta \ .\cr}
$$
Renormalizing and simplifying, we obtain the associated Askey-Wilson
recurrence
$$
X_{n+1}^{(1)}  - (z - a'_n ) X_n^{(1)}  + b_n^{\prime 2}  X_{n-1}^{(1)}  \ = \ 0
\leqno(2.5)
$$
with
$$
\eqalign{
z \ = \ & {u+ u^{-1}\over 2 } \ , \cr
a'_n \ = \ & -A'_n - B'_n + {\alpha\over 2} + {1\over 2\alpha} 
 \cr
b^{\prime 2}_n \ = \ & A'_{n-1}  B'_n  
 \cr
A'_n \ = \ &
{ (1-sq^{n-1} \epsilon ) 
  (1-\alpha\beta \epsilon q^n ) 
  (1-\alpha\gamma \epsilon q^n ) 
  (1-\alpha\delta \epsilon q^n ) \over
  2\alpha (1-sq^{2n-1} \epsilon^2) (1-sq^{2n}\epsilon^2) } \ , \cr
B'_n \ = \ &
{ \alpha (1-\epsilon q^n ) 
  (1-\beta \gamma \epsilon q^{n-1} ) 
  (1-\beta\delta \epsilon q^{n-1} ) 
  (1-\gamma\delta \epsilon q^{n-1} )  \over
  2 (1-s\epsilon^2q^{2n-2} ) (1-s\epsilon^2q^{2n-1}) }\  \cr}
$$
and a solution of (2.5) is therefore
$$
\leqalignno{
X_n^{(1)} (u) \ = \ &
 \Big( {u\over 2} \Big)^n\ 
{ \big(su\epsilon q^n/\delta\big)_\infty
 \big(s\epsilon^2q^{2n-1}\big)_\infty
\over
 \big(s\epsilon q^{n-1}\big)_\infty
 \big(\delta\epsilon q^n/u\big)_\infty
\big(\alpha \beta\epsilon q^n , 
\alpha \gamma \epsilon q^n , 
\beta \gamma \epsilon q^n \big)_\infty }& (2.6)\cr
& \ \times W \Big( {\alpha\beta \gamma u\over q}\ ; \ 
{q^{-n}\over \epsilon} , \epsilon sq^{n-1} , \alpha u, \beta u, 
\gamma u \Big) \ , \cr}
$$
with convergence condition   
$ | q/\delta u| < 1 $.

The solution  $ X_n^{(1)} $  is a generalization of the 
Askey-Wilson solution (1.1) to the  $\epsilon = 1 $  case of
(2.5).  This can be verified through an application of Watson's
formula ([5], III.18, p.242) connecting a terminating $\fphit $  with
an $ \ephis $.

A second linearly independent solution of (2.5) may be obtained with 
the help of the reflection symmetry transformation 
$ v \to -v -1 $,  
$ (\alpha, \beta, \gamma, \delta) \to 
( q /\alpha,
q /\beta,
q /\gamma,
q /\delta)$.
Using (see also [10])  
$$
\eqalign{
b_{-v-1}^2 
( q /\alpha,
q /\beta,
q /\gamma,
q /\delta)\ = \ & 
b_{ v+1}^2 (\alpha, \beta, \gamma, \delta) \cr
a_{-v-1}
( q /\alpha,
q /\beta,
q /\gamma,
q /\delta)\ = \ & 
a_v (\alpha, \beta, \gamma, \delta) \cr}
$$
and renormalizing, we arrive at the solution
$$
\leqalignno{
X_n^{(2)} (u) \ = \ &
\Big({ u\over 2}\Big)^n\ 
{ \big(s\epsilon^2 q^{2n-1}\big)_\infty
 \big(\epsilon \delta u q^{n+1}\big)_\infty
\over
 \big(\epsilon q^{n+1}\big)_\infty
 \big(\beta\delta\epsilon q^n ,
\gamma \delta\epsilon q^n , 
\alpha \delta\epsilon q^n \big)_\infty 
\big( {s\epsilon\over \delta u} q^{n-1} \big)_\infty } & (2.7)\cr
& \ \times W \Big( {q^2u \over \alpha \beta \gamma }\ ; \ 
 \epsilon q^{n+1} ,{q^{-n+2}\over \epsilon s } , {qu \over\alpha } ,
 {qu \over\beta} , 
{qu \over\gamma}  \Big) \ , \cr}
$$
with convergence condition  $ |\delta /u|< 1 $.

In order to obtain a third solution of (2.5) we make the
following alterations in (2.2):
\noindent
Interchange  $ c $  and  $ h $,  replace  $ a,b,c,\ldots , g $
by  $ h^2/ a $,  $ hb /a $,
$ hc/ a , \ldots , hg /a $  respectively, put
$ h = q^{-n} $,  $ n $  being a positive integer, and reverse
the  $ \tphin $  series.  This gives
$$
\leqalignno{
& { {b\over a} (1-h)(1-{h\over a}) 
(1-{hq\over a})
(1-{aq\over bd})
(1-{aq\over be})
(1-{aq\over bf}) 
(1-{aq\over bg}) 
(1-{aq\over bc})\over
(1- {aq\over b}) } \ 
 & (2.8) \cr
&\hskip  1in \times \big[c_2 \phi_+ (b -) - c_1 \phi\big]\cr
- & {( 1- {bh\over a} )
(1-{h\over b}) 
(1-{hq\over b})
(1-{q\over d})
(1-{q\over e})
(1-{q\over f}) 
(1-{q\over g}) 
(1-{q\over c})\over
(1- {bq\over a}) } \cr 
&\hskip  1in \times \big[c_3 \phi_- (b +) - c_1 \phi\big]\cr
- &{q\over a}
( 1- {a\over b} )
(1-{b\over q}) 
(1-{hd\over a})
(1-{he\over a})
(1-{hf\over a})
(1-{hg\over a}) 
(1-{hc\over a}) 
c_1 \phi \ = \ 0 \ , \cr}
$$
where
$$
\leqalignno{
c_1 \ = \ &
{ (\sqrt a )_n (-\sqrt a)_n 
\big( {aq\over b}\big)_n
\big( {aq\over c}\big)_n \ldots
\big( {aq\over g}\big)_n
\big( aq^{n+1} )_n\over
(a)_n (q\sqrt a)_n (-q\sqrt a)_n (b)_n (c)_n \ldots (g)_n } \ 
(-1)^n q^{n(n-1)/2} \cr
c_2 \ = \ &
{ (q\sqrt a )_{n-1} (-q\sqrt a)_{n-1} 
\big( {aq^3\over b}\big)_{n-1}
\big( {aq^2\over c}\big)_{n-1} \ldots
\big( {aq^2\over g}\big)_{n-1}
\big( aq^{n+2}\big )_{n-1}\over
(aq^2)_{n-1} (q^2\sqrt a)_{n-1} (-q^2 \sqrt a)_{n-1} (b)_{n-1} 
(cq)_{n-1} \ldots (gq)_{n-1} } \ 
(-1)^{n-1} q^{(n-1)(n-2)/2} \cr
c_3 \ = \ &
{ \big({\sqrt a \over q}\big)_{n+1} \big(-{\sqrt a\over q}\big)_{n+1} 
\big( {a\over bq}\big)_{n+1}
\big( {a\over c}\big)_{n+1} \ldots
\big( {a\over g}\big)_{n+1}
\big( aq^n \big)_{n+1}\over
\big({a\over q^2}\big)_{n+1} (\sqrt a)_{n+1} (-\sqrt a)_{n+1}
(b)_{n+1} 
\big({c\over q}\big )_{n+1} \ldots \big({g\over q}\big)_{n+1} } \ 
(-1)^{n+1} q^{n(n+1)/2} . \cr}
$$
 Note that $\phi_+(b-) $ and $\phi_-(b+) $ refer to $ \phi $ with $ (a, c, d, 
e, f, g, h) $ replaced by $(aq^2, cq, dq, eq, fq, gq, hq) $ and $(aq^
{-2}, c/q, d/q, e/q, f/q, g/q, h/q) $ respectively.

Taking the limit  $ h = q^{-n}\to \infty $,
$ g = a^3 q^{2+n} /bcdef \to 0 $,  we obtain
$$
\leqalignno{
& 
{b(1-{aq\over bd})
(1-{aq\over be})
(1-{aq\over bf})
(1-{aq\over bc})\over
a(1-{aq\over b})} & (2.9)\cr
\times & \left[
{a^2 q^2\over bcdef}\ 
{(1-aq) (1-aq^2) (1-c) (1-d) (1-e) (1-f) \over
(1-{aq\over b})
(1-{aq^2\over b})
(1-{aq\over c})
(1-{aq\over d})
(1-{aq\over e})
(1-{aq\over f})} \ 
W_+ (b-) - W \right ]\cr
& \hskip .5truein -{
(1-{q\over d})
(1-{q\over e})
(1-{q\over f}) 
(1-{q\over c})\over
(1- {bq\over a})} \cr 
\times &  \left[ {bcdef\over a^2 q^2} 
{ (1-{a\over bq})
(1-{a\over b})
(1-{a\over c})
(1-{a\over d})
(1-{a\over e})
(1-{a\over f})\over 
(1- {a\over q})  (1-a) 
(1-{c\over q})
(1-{d\over q})
(1-{e\over q})
(1-{f\over q})} \ 
W_- (b+) - W \right ]\cr
& \hskip .5truein -{q\over a} 
(1-{a\over b})
(1-{b\over q})
\big(1-{a^2 q^2\over bcdef}\big) W \ = \ 0 . \cr}
$$

In (2.9), writing
$$
\leqalignno{
a \ = \ & {\beta \gamma \delta \epsilon^2\over u} q^{2n} , \quad
b \ = \  {q\over \alpha u} , \quad
c \ = \  q^{n+1}\epsilon  , \quad
d \ = \  \beta \delta \epsilon q^n  & (2.10) \cr
e \ = \ &  \gamma \delta \epsilon  q^n, \quad
f \ = \  \beta \gamma \epsilon q^n  , \quad
s \ = \  \alpha \beta \gamma \delta \cr}
$$
and renormalizing, we obtain a third solution to (2.5) given by
$$
\leqalignno{
& \quad X_n^{(3)} (u)\ = \ 
{1\over (2u)^n} 
\cr 
& \qquad \qquad\times {\big(s\epsilon^2 q^{2n}\big)_\infty
 \big(s\epsilon^2 q^{2n-1}\big)_\infty
\big( {\beta\gamma \delta\epsilon\over u} q^n \big)_\infty 
\big( {\epsilon \beta q^{n+1}\over u},
 {\epsilon \gamma q^{n+1}\over u},
 {\epsilon \delta q^{n+1}\over u}\big)_\infty
 \over
 \big(\epsilon q^{n+1}\big)_\infty
 \big(\epsilon sq^{n-1}\big)_\infty
\big( \alpha\beta \epsilon q^n , 
\alpha \gamma\epsilon q^n,
\alpha \delta\epsilon q^n,
\beta \gamma\epsilon q^n,
\beta \delta\epsilon q^n,
\gamma \delta\epsilon q^n \big)_\infty 
\big( {\beta \gamma \delta \epsilon^2\over u} q^{2n+1} \big)_\infty }
& (2.11) \cr
&\qquad \qquad \qquad  \times W \Big( {\beta \gamma \delta \epsilon^2\over u}
q^{2n} \ ; \ 
{q\over \alpha u} , 
\epsilon q^{n+1} , 
\beta \delta\epsilon q^n , 
\gamma \delta \epsilon q^n ,
\beta \gamma \epsilon q^n  \Big) \ , \cr}
$$
with convergence condition  $ | \alpha/ u | < 1 $.
A transformation ([5], III.23, p.243) may be applied to  $ X_n^{(3)}$
to obtain this solution in a form which is symmetric in the
parameters  $ \alpha, \beta, \gamma, \delta $.  We get the solution
(omitting constant factors)
$$
\leqalignno{
& X_n^{(4)} (u) \cr
= \ &
{1\over (2u)^n}\ 
{\big(\epsilon ^2 sq^{2n-1}\big)_\infty
\big( {\alpha\epsilon\over u} q^{n+1} ,
 {\beta \epsilon \over u}  q^{n+1},
 {\gamma \epsilon \over u}  q^{n+1},
 {\delta \epsilon \over u}  q^{n+1}\big)_\infty
 \over
 \big(\epsilon q^{n+1}\big)_\infty
 \big({q^{n+2}\epsilon\over u^2}\big)_\infty
\big( \alpha\beta \epsilon q^n , 
\alpha \delta\epsilon q^n,
\alpha \gamma\epsilon q^n,
\beta \gamma\epsilon q^n,
\beta \delta\epsilon q^n,
\gamma \delta\epsilon q^n \big)_\infty }
& (2.12) \cr
&\ \times W \Big( 
{q^{n+1}\epsilon \over u^2} ; \ 
\epsilon q^{n+1},
{q\over \alpha u} , 
{q\over \beta u},
{q\over \gamma u},
{q\over \delta u} \Big)  \cr}
$$
where the convergence condition is  $ |sq^{n-1}\epsilon | < 1 $.

The large  $ n $ behaviour of  $ X_n^{(4)} (u) $  is easily seen
to be
$$
\leqalignno{
& X_n^{(4)} (u)\ 
\ttopsim{n\to\infty}
{1\over (2u)^n}\cr
\noalign{\hbox{while the corresponding solution}}
& X_n^{(4)} (1/u)
\ \sim \ \big({u\over 2}\big)^n . \cr}
$$
It is evident that  $ X_n^{(4)} (u) $  is a sub-dominant or
minimal solution
for  $ |u| > 1 $,  $ |sq^{n-1} \epsilon | < 1 $   while 
$ X_n^{(4)} (1 /u) $ 
is sub-dominant (minimal) for $ | u | < 1 $,
$ |sq^{n-1} \epsilon | < 1 $.

We proceed to get another solution by making the following parameter
replacements in (2.9):
$$
\eqalign{
a \ = \ & {q^{-2n+1}\over \beta\gamma \delta u\epsilon^2}, \quad
b \ = \  {\alpha \over u},\quad
c \ = \  {q^{-n}\over \epsilon } ,\quad
d \ = \  {1\over \beta\delta \epsilon } \ q^{-{n+1}},\cr
e \ = \  & {1\over \gamma\delta \epsilon } \ q^{-{n+1}},\quad
f \ = \   {1\over \beta \gamma \epsilon } \ q^{-{n+1}}\ .\cr}
\leqno(2.14)
$$
Note that this is just the reflection transformation that was 
used to obtain the solution $ X_n^{(2)} $ from $ X_n^{(1)} $ but
now it is being applied to $ X_n^{(3)} $.

After renormalization, omitting constant factors
and using the transformation [5, III.23, p.243] we arrive at the
following solution 
$$
\leqalignno{
 X_n^{(5)} (u) 
\ = \ & {1\over (2u)^n}\ 
{\big(sq^{2n-1}\epsilon^2\big)_\infty
\big( u^2 \epsilon q^n\big)_\infty
\over 
\big(sq^{n-1}\epsilon\big)_\infty
\big(\alpha u \epsilon q^n , 
\beta u\epsilon q^n,
\gamma u \epsilon q^n,
\delta u \epsilon q^n\big)_\infty } & (2.15)\cr
&\quad \times W \Big( 
{q^{-n}\over \epsilon  u^2} ; \ 
{q^{-n}\over \epsilon },
{\alpha \over u} , 
{\beta\over u},
{\gamma\over u},
{\delta \over u} \Big)  \cr}
$$
where the convergence condition is
$ \big| q^{-n+2} /s\epsilon \big| < 1 $.

Note that, when $ \epsilon  =1 $, the solution (2.15) is proportional
to the Askey-Wilson polynomial $ P_n(z; \alpha ,\beta , \gamma , \delta
) $ in (1.1) and can be used to give $ P_n(z) $ in a form which is
explicitly symmetric in its parameters. That is
$$
P_n(z) \ = \ (2u)^{-n}{(\alpha u, \beta u, \gamma u, \delta u, s/q)_n
\over (s/q)_{2n}(u^2)_n}W(q^{-n}/u^2;q^{-n},{\alpha \over u},
{\beta \over u}, {\gamma \over u}, {\delta \over u}) .
$$

A solution in which the convergence condition is independent
of  $ z $  and also independent of the parameters
$ \alpha, \beta, \gamma, \delta $  can be obtained from (2.12)
by applying the transformation
([5], III.24, p.243).  This gives
$$
\leqalignno{
 X_n^{(6)} (u) 
\ = \ & {1\over (2u)^n}\ 
{\big(s\epsilon^2 q^{2n-1}\big)_\infty \big(
s\epsilon q^n/ \alpha u,
s\epsilon q^n/ \beta  u,
s\epsilon q^n/ \gamma  u,
s\epsilon q^n/ \delta  u\big)_\infty\over
\big(\epsilon sq^{n-1}\big)_\infty
\big(s\epsilon  q^n/u^2 \big)_\infty
\big( \alpha \beta\epsilon q^n , 
\alpha\gamma \epsilon q^n,
\alpha\delta \epsilon q^n,
\beta\gamma \epsilon q^n,
\beta\delta \epsilon q^n,
\gamma\delta \epsilon q^n\big)_\infty} \cr
&\quad \times W \Big( 
{s\epsilon q^{n-1}\over u^2} ; \ 
{\epsilon s q^{n-1}},
{\alpha \over u} , 
{\beta\over  u},
{\gamma\over  u},
{\delta \over u} \Big) & (2.16) \cr}
$$
with convergence condition  $ \big| q^{n+1}\epsilon \big| < 1 $.

Since  $ X_n^{(6)} (u) \topsim{n\to\infty} 
{1\over (2u)^n} $,
 we have that $  X_n^{(6)} (u) $  is a sub-dominant(minimal) solution for
$ |u | > 1 $  while
  $  X_n^{(6)} (1/u) $  is a sub-dominant(minimal) solution for 
$ |u | < 1 $.
 Note that  $ X_n^{(6)} $  is really the same solution as
$ X_n^{(4)} $,  but analytically continued.

Summarizing we have

\proclaim Theorem 1.  The associated Askey-Wilson equation (2.5)
has solutions  $ X_n^{(i)} (u) $ and $ X_n^{(i)}(1/u) $ ,
$ i = 1,2,\ldots , 6 $  given by
(2.6), (2.7), (2.11), (2.12), (2.15) and (2.16).  If  $ |u| > 1 $,  
then  
$ X_n^{(3)} (u) $,   
$ X_n^{(4)} (u) $  and  
$ X_n^{(6)} (u) $ 
(which are connected by  general  $ \ephis $  transformations) each
represent a minimal solution
$ X_n^{(s)} (z) $ 
for the parameter values  
$ \big|{\alpha\over u}\big| < 1 $,
$ |\epsilon sq^{n-1}| < 1 $,
$ |\epsilon q^{n+1}| < 1 $
respectively.\par
\vglue .5truein

\beginsection 3. Continued fraction representations.

The continued fraction associated with the difference equation
(2.5) is
$$
CF (z) \ :=\ z - a'_0 + \bigk\limits_{n=1}^\infty
\left( {-b_n^{\prime 2} \over z -a'_n } \right)  \ .
\leqno(3.1)
$$
If  $ b_n^{\prime 2} \not= 0 $,  $ n \ge 1 $,  by Pincherle's
theorem [6], [17] and Theorem 1 of Section 2 we have
$$
{1\over CF (z)} \ =\ 
{ X_0^{(s)} (z)\over
b_0^{\prime 2} X_{-1}^{(s)} (z)} \ .
\leqno(3.2)
$$
Using  
$ X_0^{(s)} (z) = X_0^{(4)} (u)$,
$ X_{-1}^{(s)} (z) = X_{-1}^{(4)} (u)$,
$ |u| > 1, z={ u + u^{-1} \over 2} $,
we have from (2.12),
$$
\leqalignno{ 
{1\over CF (z)} \ =\ 
& {2\over u} 
{ (1-\epsilon q/ u^2) 
(1- s\epsilon^2/ q^2)
(1- s\epsilon^2/ q)
\over
( 1-\alpha\epsilon / u)
( 1-\beta\epsilon/ u)
( 1-\gamma\epsilon/ u)
( 1-\delta\epsilon/ u)
( 1-s\epsilon/ q^2) }\cr
& \times 
{W \Big( 
q\epsilon/ u^2 ; \ q\epsilon,
q/ \alpha u,
q/ \beta u,
q/ \gamma  u,
q/ \delta  u\Big)
\over
W \Big( 
\epsilon/ u^2 ; \ \epsilon,
q/ \alpha u,
q/ \beta u,
q/ \gamma  u,
q/ \delta  u\Big)} & (3.3)\cr}
$$
for  $ |u| > 1 $  and  $ |s\epsilon/ q^2 | < 1 $.
For $|u| < 1, \,|s\epsilon/q^2| < 1 $ we replace $u$ in
(3.3) by $ 1/u $.

 Other representations for the continued fraction for 
different parameter ranges are
obtained by taking a different representation for the minimal
solution.
 For example, by taking  
$ X_0^{(s)} (z) = X_0^{(6)} (u)$ and
$ X_{-1}^{(s)} (z) = X_{-1}^{(6)} (u)$, 
$ | u | > 1 $,  we obtain the continued fraction representation
$$
\leqalignno{ 
{1\over CF (z)} \ =\ 
& {2\over u}{ 
(1- s\epsilon^2 /q^2)
(1- s\epsilon^2 /q)(1-s\epsilon/u^2q) \over
(1-\epsilon)(1-s\epsilon/\alpha uq)(1-s\epsilon/\beta uq)(1-s\epsilon
/\gamma uq)(1-s\epsilon/\delta uq)}& (3.4)\cr
& \times
{W \Big( 
\epsilon s/u^2q  ; \ \epsilon s/q,
 \alpha /u,
 \beta/ u,
 \gamma / u,
 \delta / u\Big)
\over
W \Big( 
s\epsilon /u^2q^2 ; \ \epsilon s /q^2 ,
\alpha /u,
\beta /u,
\gamma / u,
\delta / u\Big)} , \cr}
$$
for  $ |u| > 1 $,  $ | \epsilon |  < 1 $.
\vglue .5truein

\beginsection 4. Weight function.

We have seen above that a minimal solution exists for
$ | u | > 1 $  and  $ | u | < 1 $.  There is no minimal solution
for  $ z \in [-1,1] $,  $(|u|=1)$,  where we have a continuous
spectrum for the associated tridiagonal Jacobi matrix  $ J $
with diagonal 
$ (a'_0, a'_1, \ldots) $  and  
$ (b'_1, b'_2, \ldots) $  above and below the diagonal.
For a probability measure  $ d\omega (x;\epsilon) $
we have, in the case of real orthogonality $ ( b_{n+1}^{\prime 2} 
> 0 , a'_n\quad\hbox{real},\quad n\ge 0 )$ for the associated monic
$ q $-Askey-Wilson polynomials
$ P_n (x; \epsilon)$,
$$
\int_\R \ P_n (x;\epsilon) P_m (x;\epsilon)
d\omega (x;\epsilon) \ = \ 
\delta_{nm}  \prod\limits_{k=1}^n \ b_k^{\prime 2} \ .
\leqno(4.1)
$$
From [21], we have the representation
$$
{1\over CF(z)} \ = \ 
\int_\R \ 
{d\omega (x;\epsilon)\over z-x} \ = \ 
{X_0^{(s)} (z)\over b_0^{\prime 2} X_{-1}^{(s)} (z)}
\leqno(4.2)
$$
and for the absolutely continuous part for  $ x\in [-1,1] $,
$$
\leqalignno{
{d\omega (x;\epsilon)\over dx} \ = \ 
& {1\over 2\pi i b_0^{\prime 2} } \ 
\left( 
{X_0^{(s)} (x-i 0)\over 
X_{-1}^{(s)} (x-i 0)} -
{X_0^{(s)} (x+i 0)\over 
X_{-1}^{(s)} (x+i 0)} \right) & (4.3) \cr
= \ & {1\over 2\pi i  b_0^{\prime 2} } \ 
{\cW \big(X_{-1}^{(s)} (x+i 0),
X_{-1}^{(s)} (x-i 0) \big)\over
\big| X_{-1}^{(s)} (x-i 0)\big|^2 } , \cr}
$$
where  
$$
W (X_n, Y_n ) := X_n Y_{n+1} - X_{n+1} Y_n \ .
$$

From (2.5) we have
$$
\cW \Big(  X_{-1}^{(4)} (u), X_{-1}^{(4)} 
(1/ u)\Big) \ = \ 
\lim\limits_{n\to \infty} 
{\cW \Big( X_n^{(4)} (u), X_n^{(4)} 
(1 /u)\Big) \over
\prod\limits_{k=0}^n b_k^{\prime 2}}.
\leqno(4.4)
$$
Using (2.5) and (2.12) we then obtain
$$
\leqalignno{
\cW \Big(  X_{-1}^{(4)} (u), & X_{-1}^{(4)} 
(1 /u)\Big) & (4.5) \cr
& = \  2(u-u^{-1})  
{\big({s\epsilon^2\over q^3}\big)_\infty
\big({s\epsilon^2\over q^2}\big)_\infty
\over
\Big( {\alpha \beta \epsilon \over q} ,
{ \alpha \gamma \epsilon \over q} ,
{ \alpha \delta \epsilon \over q} ,
{ \beta \gamma \epsilon \over q} ,
{ \beta \delta \epsilon \over q} ,
{ \gamma \delta\epsilon\over q} ,
{ s \epsilon\over q^2} ,
\epsilon \Big)_\infty }\ . \cr}
$$
Thus (4.3) with $$ X_n^{(s)}(x \pm i0)\ = \ X_n^{(4)}(e^{\mp i\theta}),
\quad
x\ = \ \cos\theta,\quad u \ = \ e^{i\theta} $$

 \noindent gives for  $ -1\le x \le 1 $,  
$ \big| s\epsilon /q^2 \big| < 1 $,
$$
\leqalignno{
& {d\omega (x;\epsilon )\over dx} \ = \ 
{2\sqrt {1-x^2}\over \pi }\ 
{(1 - s\epsilon^2/q)
(1 - s\epsilon^2 /q^2)^2\over
(1-s\epsilon/q^2)} & (4.6)\cr
\times \ &
{\big(
\epsilon q /u^2, \epsilon q u^2,
\alpha \beta \epsilon,
\alpha \gamma \epsilon,
\alpha \delta \epsilon,
\beta\gamma \epsilon,
\beta\delta \epsilon,
\gamma\delta \epsilon,
\epsilon q \big)_\infty\over
\big(\alpha\epsilon /u ,
\alpha\epsilon u,
\beta \epsilon /u,
\beta \epsilon u,
\gamma\epsilon /u,
\gamma\epsilon u,
\delta\epsilon /u,
\delta\epsilon u,
s\epsilon /q^2\big)_\infty}  \cr
\times \ &
{1\over \Big| W \big( 
\epsilon /u^2; \ 
q /\alpha u,
q /\beta u,
q /\gamma u,
q /\delta u,
\epsilon \big)\Big|^2} \ .\cr}
$$
This checks with the weight function obtained by Ismail and
Rahman ([16], (4.31), p.218). When  $ \epsilon = 1 $,  this
reduces to the Askey-Wilson weight function [1].
Also from the positivity of the denominators in (3.4) there
is no discrete spectrum if   $ -1 < \epsilon < 1$,
$ |\alpha |, | \beta|, |\gamma |, |\delta |< |q|^{1/2} $, [16].
Other conditions for the absence of a discrete spectrum
may be deduced from (3.3) and (3.4). For example, from (3.3),
there is no discrete spectrum if $ -|q^2/s| < \epsilon < |q^2/s| $,
$ | \alpha |, |\beta |, |\gamma |, |\delta | > |q|^{1/2} $ .

Next we use (4.3) with
$$
{d\omega (x;\epsilon )\over dx} \ = \ 
{1\over 2\pi i b_0^{\prime 2}}
\left( 
{X_0^{(6)} (u)\over 
X_{-1}^{(4)} (u)} -
{X_0^{(6)} (1/u)\over 
X_{-1}^{(4)} (1/u)} \right) \ ,\quad x \ = \ \cos\theta,\quad u \ 
= \ e^{i\theta}.
$$
The right side simplifies to
$$
\leqalignno{
& {1\over 2\pi i}
{\big(1 - s\epsilon^2/ q^2\big)
\big(1 - s\epsilon^2/ q\big) ( \epsilon q)_\infty\over
\big(s\epsilon/ q^2\big)_\infty } & (4.7)\cr
\times \ & 
\Bigg[ 2u
{( \epsilon q u^2)_\infty
\Big( {s\epsilon\over \alpha} u ,
{s\epsilon\over \beta} u ,
{s\epsilon\over \gamma} u ,
{s\epsilon\over \delta} u \Big)_\infty
W \big({su^2 \epsilon\over q} ; \ 
{s\epsilon\over q},
\alpha u, \beta u, \gamma u, \delta u\big)\over
(su^2 \epsilon)_\infty
(\alpha u\epsilon, \beta u\epsilon, \gamma u\epsilon, 
\delta u\epsilon)_\infty
W \big( \epsilon u^2;\ \epsilon, 
{qu\over \alpha},
{qu\over \beta},
{qu\over \gamma},
{qu\over \delta}\big)} \cr
&  - {2\over u}
{\big(\epsilon q/ u^2\big)_\infty
\big( s\epsilon/ \alpha u,
s\epsilon /\beta u,
s\epsilon /\gamma u,
s\epsilon /\delta u\big)_\infty
W \big( {s\over u^2}
{\epsilon\over q}; \ 
{s\epsilon\over q}, \ 
{\alpha\over u},
{\beta \over u},
{\gamma \over u},
{\delta \over u}\big)\over
\big( {s\epsilon\over u^2}\big)_\infty
\big( {\alpha\epsilon\over u},
{\beta \epsilon\over u},
{\gamma\epsilon\over u},
{\delta \epsilon\over u}\big)_\infty
W \big( {\epsilon\over  u^2};\ \epsilon, 
{q\over \alpha u},
{q\over \beta u},
{q\over \gamma u},
{q\over \delta u}\big)} \Bigg] . \cr}
$$
Equating (4.6) and (4.7), and writing
$$
\leqalignno{
& G (\alpha, \beta, \gamma , \delta, \epsilon, u) := {1\over u}\ 
{\big( {\epsilon q\over u^2}\big)_\infty
\big( 
{s\epsilon\over \alpha u},
{s\epsilon\over \beta  u},
{s\epsilon\over \gamma  u},
{s\epsilon\over \delta u}\big)_\infty\over
\big(s\epsilon/ u^2\big)_\infty
\big(
{\alpha\epsilon\over u},
{\beta\epsilon\over u},
{\gamma \epsilon\over u},
{\delta \epsilon\over u}\big)_\infty} & (4.8)\cr
\times & 
W \big( {s\epsilon\over u^2 q} ; \ 
{s\epsilon\over q},
{\alpha\over u},
{\beta \over u},
{\gamma \over u},
{\delta \over u}\big)
W \big( \epsilon u^2;\ \epsilon, 
{qu\over \alpha},
{qu\over \beta},
{qu\over \gamma},
{qu\over \delta}\big)\ , \cr}
$$
we obtain the identity
$$
\leqalignno{
& G (\alpha, \beta, \gamma, \delta , \epsilon, u) - 
 G \big(\alpha, \beta, \gamma, \delta , \epsilon, 1/ u\big ) 
 \ = \ \big(1/ u - u\big) \ 
{(1-s\epsilon^2/ q^2 )\over
(1-s\epsilon/ q^2)}
 & (4.9)\cr
\times \ & 
{\big
(\alpha \beta \epsilon,
\alpha \gamma \epsilon,
\alpha \delta \epsilon,
\beta\gamma \epsilon,
\beta\delta \epsilon,
\gamma\delta \epsilon,
\epsilon q/ u^2,
\epsilon q u^2\big)_\infty\over
\big(\alpha\epsilon/ u ,
\alpha\epsilon u,
\beta \epsilon/ u,
\beta \epsilon u,
\gamma\epsilon/ u,
\gamma\epsilon u,
\delta\epsilon/ u,
\delta\epsilon u\big)_\infty }\cr}
$$
which is a  $ q $-analogue of Masson's generalization of Dougall's
theorem [18].  $ \epsilon = 1 $  gives a  $ q $-analogue of 
Dougall's theorem (see [3]).

We can recover from (4.9) the following identity which we had 
obtained in our earlier paper ([7], (31), p.723) for
$ \epsilon = 1 $,  $ s = q^m $;  $ m = 1,2,\ldots , $ i.e.
$$
\leqalignno{
& \big({q\over \alpha}\big)^{m-3}
\big( {1\over u} - u \big)
\Big[ u^{m-2} \Pi_1 (u) \Pi_2 \big({1\over u}\big) -
u^{2-m} \Pi_1 \big({1\over u}\big) \Pi_2 (u)\Big] & (4.10)\cr
& = ( \alpha \beta q^{-1}, 
\alpha\gamma q^{-1},
\alpha \delta q^{-1})_\infty
\big( {q^2\over \alpha\beta},
 {q^2\over \alpha\gamma },
 {q^2\over \alpha\delta }\big)_\infty
 \big(u^2, {1\over u^2}\big)_\infty \cr}
$$
for  $ m = 1,2, \ldots$ and $ |u| = 1 $,  where
$$
\leqalignno{
\Pi_1(u) \ = \ & \big( 
{q\over \alpha u},
{q\over \beta  u},
{q\over \gamma  u},
{q\over \delta  u}\big)_\infty \cr
\Pi_2(u) \ = \ & \big( 
{\alpha\over u},
{\beta\over  u},
{\gamma \over u},
{\delta\over  u}\big)_\infty\ . \cr}
$$
\vfill\eject

\centerline {\bf References }
\baselineskip=12pt
\frenchspacing
\vglue .25in
\item{1. }
R. Askey and J. Wilson,
Some basic hypergeometric orthogonal polynomials that generalize
Jacobi polynomials,
{\it Memoirs Amer. Math. Soc.} {\bf 319} (1985) 1--55.
\medskip

\item{2. }
R. Askey and J. Wimp,
Associated Laguerre and Hermite polynomials,
{\it Proc. Roy. Soc. Edinburgh} Sect. A 96 (1984), 15--37.
\medskip

\item{3. }
W.N. Bailey,
{\it Generalized Hypergeometric Series,}
Cambridge Univ. Press, London, 1935.
\medskip

\item{4. }
J. Bustoz and M.E.H. Ismail,
The associated ultraspherical polynomials and their $ q $-analogues,
{\it Canad. J. Math.} {\bf 34} (1982), 718--736.
\medskip

\item{5. } G. Gasper and M. Rahman,
{\it Basic Hypergeometric Series,}
Cambridge Univ. Press, Cambridge, 1990.
\medskip

\item{6. }
W. Gautschi,
Computational aspects of three-term recurrence relations,
{\it SIAM Rev.} {\bf 9} (1967), 24--82.
\medskip

\item{7. } D.P. Gupta and D.R. Masson,
Exceptional  $ q $-Askey-Wilson polynomials and continued fractions,
{\it Proc. A.M.S.} {\bf 112} (1991), 717--727.
\medskip

\item{8. } D.P. Gupta, M.E.H. Ismail and D.R. Masson,
Associated continuous Hahn polynomials,
{\it Canad. J. of Math.} {\bf 43} (1991), 1263--1280.
\medskip

\item{9. } D.P. Gupta, M.E.H. Ismail and D.R. Masson,
Contiguous relations, Basic Hypergeometric functions and
orthogonal polynomials II, Associated big $ q $-Jacobi polynomials,
{\it J. of Math. Analysis and Applications}  {\bf 171} (1992), 477--497.
\medskip

\item{10. } 
D.P. Gupta and D.R. Masson,
Watson's basic analogue of Ramanujan's Entry 40 and its
generalization,
{\it SIAM J. Math. Anal.}, to appear.
\medskip

\item{11. } M.E.H. Ismail, J. Letessier, and G. Valent,
Linear birth and death models and associated Laguerre polynomials,
{\it J. Approx. Theory} {\bf 56} (1988), 337--348.
\medskip

\item{12. } M.E.H. Ismail, J. Letessier, and G. Valent,
Quadratic birth and death processes and associated continuous
dual Hahn polynomials,
{\it SIAM J. Math. Anal.} {\bf 20} (1989), 727--737.
\medskip

\item{13. } M.E.H. Ismail, J. Letessier, G. Valent and J. Wimp,
Two families of associated Wilson polynomials,
{\it Can. J. Math.} {\bf 42} (1990), 659--695.
\medskip

\item{14. } M.E.H. Ismail and C.A. Libis,
Contiguous relations, basic hypergeometric functions and orthogonal
polynomials I,
{\it J. Math. Anal. Appl.} {\bf 141} (1989), 349--372.
\medskip

\item{15. } M.E.H. Ismail and D.R. Masson,
Two families of orthogonal polynomials related to Jacobi
polynomials,
{\it Rocky Mountain J. Math.} {\bf 21} (1991), 359--375. 
\medskip

\item{16. } M.E.H. Ismail and M. Rahman,
Associated Askey-Wilson polynomials,
{\it Trans. Amer. Math. Soc.} {\bf 328} (1991), 201--239.
\medskip

\item{17. }
W.B. Jones and W.J. Thron, {\it Continued Fractions: Analytic
Theory and Applications,} 
Addison-Wesley, Reading, Mass., 1980.
\medskip

\item{18. } D.R. Masson, 
Associated Wilson polynomials, 
{\it Constructive Approximation} {\bf 7} (1991), 521--534.
\medskip

\item{19. }
\bibline ,
Wilson polynomials and some continued fractions
of Ramanujan,
{\it Rocky Mountain J. of Math.} {\bf 21} (1991), 489--499.
\medskip

\item{20. }
\bibline ,
The rotating harmonic oscillator eigenvalue problems, I.
Continued fractions and analytic continuation,
{\it J. Math. Phys.} {\bf 24} (1983), 2074--2088.
\medskip

\item{21. }
\bibline,
Difference equations, continued fractions, Jacobi Matrices and
orthogonal polynomials, In:
{\it Non-linear numerical methods and Rational Approximation}
(A. Cuyt, ed.) Dordrecht, Reidel, 1988, 239--257.
\medskip

\item{22. }
J.A. Wilson,
Hypergeometric series, recurrence relations and some new orthogonal
polynomials,
Ph.D. diss., University of Wisconsin, Madison, 1978.
\medskip

\item{23. }
J. Wimp,
Explicit formulas for the associated Jacobi polynomials and
some applications,
{\it Canad. J. Math.} {\bf 39} (1987), 983--1000.
\end